\documentclass[12pt]{article}
\usepackage{amsmath,amsfonts,amscd,amssymb}
\newcommand{\IR}{{\mathbb R}} % !!!!!!!!!!!!!!!!!!!!!!!!!!!!!!!!!
\font\caps=cmcsc10 scaled\magstep1
\font\sefont=cmr12 scaled \magstep1 %12 statt 17
\def\bs{\bigskip}
\def\noi{\noindent}

\def\ms{\medskip}
\def\ve{\varepsilon}
\def\vf{\varphi}

\def\Strich{\par\ms
            \centerline{\hbox to 2cm {\hrulefill}}\par\bs\noi}
\def\qed{\hfill{$\square$}}
\def\de{\partial}

\font\nbf=cmbxti10 scaled\magstep1                       %fett kursiv

\begin{document}
%%%%%%%%%%%%%Ende Vorspann von Frau Zirngibl%%%%%%%%%%

\thispagestyle{empty}
\def\dg{DGL($\lambda$)\ }
\def\eps{\ve}
\def\Tr{\mathop{\rm Tr}}
\def\tr{\mathop{\rm tr}}
\def\inv{{}^{-1}}
\def\id{\rm id}

\centerline{\bf\Large
Mathematical remarks on}

\ms
\centerline{\bf\Large
transcritical bifurcation in Hamiltonian systems}

\vskip.7cm
\centerline{K. J\"anich} %, September 4, 2005 [26.7.06]}
\vskip.5cm
\centerline{\it Faculty of Mathematics, University of Regensburg,
            D-93040 Regensburg, Germany}
\vskip0.8cm

\centerline{\bf Abstract}

\bigskip

\baselineskip 10pt

\noindent
This article is meant as a mathematical appendix or comment on
\cite{BT}. We first consider the notion of transcritical
bifurcations of fixed points of general area-preserving maps,
and then adress some questions related to \cite{BT} on bifurcation 
in Poincar\'e maps of 2-dimensional Hamiltonian systems.

\baselineskip 14pt

\vskip1.cm
\noi
{\sefont 1 Rank-1-bifurcations}

\ms\noi
Let differentiable ($C^\infty$) functions
$Q=Q(q,p,\eps)$ and
$P=P(q,p,\eps)$
be defined on an open neighborhood of the origin in $(q,p,\eps)$-space
$\IR^3$. We assume $(P,Q)$ to be symplectic in the
$(q,p)$-coordinates, which means that
\begin{equation}
\det\left(
\begin{array}{ll}
Q_q(q,p,\eps)&Q_p(q,p,\eps)\\
P_q(q,p,\eps)&P_p(q,p,\eps)\\
\end{array}
\right)
\equiv
1,
\end{equation}
where we have written the partial derivatives as
$\frac{\de Q}{\de q}=:Q_q$ etc. We then speak of $(Q,P)$ as
describing a {\nbf symplectic family}.
Our bifurcating fixed point shall be the origin of
the $(q,p)$-plane at the value $\eps=0$ of the parameter, so
$Q(0,0,0)=0$ and $P(0,0,0)=0$. At the bifurcation point
$(0,0,0)$, we assume the Jacobian matrix with respect to the
variables $(q,p)$ to have the eigenvalue $+1$. Other\-wise, by the
implicit function theorem, the fixed point set could be
parametrized locally by $\eps$ and hence would not bifurcate in
the way we wish to study. We will exclude, however, the
exceptional case that the Jacobian matrix equals the identity
matrix. So by assumption, the eigenspace is a 1-dimensional
subspace of the $(q,p)$-plane. We can
always adjust the canonical coordinates, for example by a simple rotation,
to have the eigenspace as the $q$-axis. Then at $(0,0,0)$ we
have
\begin{equation}
\left(
\begin{array}{ll}
Q_q&Q_p\\
P_q&P_p\\
\end{array}
\right)
=
\left(
\begin{array}{lc}
1&Q_p\\
0&1\\
\end{array}
\right)\quad\text{and}\quad Q_p\ne0
\end{equation}
in these coordinates. The total
fixed point set
\begin{equation}
F:=\{(q,p,\eps)\mid Q(q,p,\eps)=q;\;P(q,p,\eps)=p\}
\end{equation}
is the inverse image of the origin $(0,0)$ in $\IR^2$ under the map
$(Q-q,P-p)$, and the Jacobian matrix of this map at $(0,0,0)$ is
\begin{equation}
\left(
\begin{array}{ccc}
Q_q\mathord-1&Q_p&Q_\eps\\
P_q&P_p\mathord-1&P_\eps\\
\end{array}
\right)
=
\left(
\begin{array}{ccc}
0&Q_p&Q_\eps\\
0&0&P_\eps\\
\end{array}
\right)\,.
\end{equation}
If we were to look for generic bifurcations, we could now impose
the condition that this matrix be of rank 2, or $P_\eps\ne0$ at
$(0,0,0)$. Then $F$ would be locally a smooth 1-dimensional
submanifold, tangent to the $q$-axis at the bifurcation point,
and if we require, as a further generic condition, $\eps|F$ to
have a nondegenerate extremum at this point, we would arrive at the
notion of an {\it extremal fixed point} or
{\it saddle-node bifurcation} in the sense of K.~R.~Meyer \cite{M}.

Instead, in this note we will be interested in bifurcations with
matrix (4) being of rank~1.
Among these {\nbf rank-1-bifurcations} we will distinguish
various types by further conditions. For a convenient coordinate
description, let us adjust the canonical $(q,p)$-coordinates
one step further by an $\eps$-dependent translation
\begin{equation}
\begin{array}{lll}
\widetilde q&=&q\\
\widetilde p&=&p-c\eps\\
\end{array}
\end{equation}
where the constant $c$ is defined by $c:=-Q_\eps/Q_p$ at $(0,0,0)$.
Then the
$\widetilde q$-axis still is the eigenspace at the bifurcation
point, and the Jacobian matrix of
$(\widetilde Q-\widetilde q,\widetilde P-\widetilde p)$ at
$(0,0,0)$ is simplified to
\begin{equation}
\left(
\begin{array}{ccc}
\widetilde Q_{\widetilde q}\mathord-1&\widetilde Q_{\widetilde
p}&\widetilde Q_\eps\\
\widetilde P_{\widetilde q}&\widetilde P_{\widetilde
p}\mathord-1&\widetilde P_\eps\\
\end{array}
\right)
=
\left(
\begin{array}{ccc}
0&Q_p&0\\
\noalign{\vskip5pt}
0&0&0\\
\end{array}
\right)\,.
\end{equation}
So let $(q,p)$ be such coordinates to begin with. To have a name
for them, we introduce the following terminology.

\bs\noi
{\bf Definition 1:} Symplectic coordinates $(q,p)$ shall be
called {\nbf adapted coordinates} for a rank-1-bifurcation at
$(0,0,0)$, if
\begin{equation}
\left(
\begin{array}{ccc}
 Q_{ q}\mathord-1& Q_{
p}& Q_\eps\\
 P_{ q}& P_{
p}\mathord-1& P_\eps\\
\end{array}
\right)
=
\left(
\begin{array}{ccc}
0&Q_p&0\\
\noalign{\vskip5pt}
0&0&0\\
\end{array}
\right)
\end{equation}
at the bifurcation point, and $Q_p\ne0$ there.\qed

\bs\noi
Then the
inverse image of 0 under $Q-q$, let's denote it by
\begin{equation}
X:=\{(q,p,\eps)\mid Q(q,p,\eps)=q\},
\end{equation}
will locally at $0=(0,0,0)$ be a smooth surface, with the
$(q,\eps)$-plane as tangent plane $T_0X$ at this point. The fixed point
set $F$ is contained in $X$, it is the inverse image
\begin{equation}
F=\{(q,p,\eps)\in X\mid P(q,p,\eps)=p\},
\end{equation}
of 0 under the restriction $(P-p)|X$ of $P-p$ to $X$. As
$(0,0,0)$ is a critical point of this restriction, we naturally
turn to its Hessian quadratic form on the $(q,\eps)$-plane
for information about the
local behavior of $P-p$ on $X$. And since $(0,0,0)$ is
critical not only for the restriction, but also for the
function $P-p$ itself, we know that this quadratic form is
simply described by the Hessian matrix of $P-p$, hence of $P$,
at $(0,0,0)$ with respect to the $(q,\eps)$-coordinates. This is
how this matrix enters the following definition.

\bs\noi
{\bf Definition 2:} A rank-1-bifurcation shall be called
{\nbf regular}, if in adapted coordinates the Hessian
\begin{equation}
\left(
\begin{array}{ll}
P_{qq}&P_{q\eps}\\
P_{q\eps}&P_{\eps\eps}\\
\end{array}
\right)\quad\text{at}\quad(0,0,0)
\end{equation}
of $P$ with respect to $q$ and $\eps$ is nondegenerate.
Depending on whether it is definite or indefinite, we speak of the
rank-1-bifurcation as being {\nbf definite} or
{\nbf indefinite}. An
indefinite rank-1-bifurcation
will be called a {\nbf cross-bifurcation},
and a cross-bifurcation
is called {\nbf transcritical}, if
the $q$-axis is not contained in the zero set of the quadratic
form, that is if $P_{qq}\ne0$ at $(0,0,0)$, in adapted
coordinates. \qed

\bs\noi
The definition does not depend on the {\it choice} of the adapted
coordinates. As can be shown, the Hessians from different
choices of adapted coordinates are equivalent (up to sign) as
quadratic forms.

%\break
\vskip1cm
\noi
{\sefont 2 The fixed point set: branches and traces}

\ms\noi
Let a regular rank-1-bifurcation in adapted coordinates be
given. We now denote the restriction $(P-p)|X$ by
$\psi:X\to\IR$. So we are interested in $\psi\inv(0)\subset X$,
because this is the fixed point set $F$ of the bifurcation.
The Hessian quadratic form
\begin{equation}
{\rm Hess}_\psi:T_0X\to\IR
\end{equation}
of $\psi$ at the bifurcation point is given by the matrix (10).
In coordinates $(q,\eps)$ on $X$,
Taylor expansion to second order of $\psi$ at this point gives
$\psi\approx\frac12{\rm Hess}_\psi$ up to higher order terms. But since
the Hessian is assumed to be nondegenerate, we can do better
than that. By the Morse Lemma, see for instance \cite{PS}, we can
find local coordinates on $X$, approximating $(q,\eps)$ in first
order, in which $\psi$ actually {\it coincides} with
$\frac12{\rm Hess}_\psi$. More precisely, there is a
diffeomorphism $f:\Omega\to \Omega'$ from an open
neighborhood $\Omega$ of the
origin in the tangent plane of $X$ to an open neighborhood $\Omega'$
of the bifurcation point in $X$ itself, such that $f(0)=0$ and
\begin{equation}
\psi(f(v))=\textstyle\frac12{\rm Hess}_\psi(v)
\end{equation}
for all $v\in \Omega$, and the differential $df_0:T_0X\to T_0X$ is the
identity. Thus up to this diffeomorphism, the fixed point set of
the bifurcation looks locally the same as the zero set
${\rm Hess}_\psi^{-1}(0)$ of the Hessian, which is a single point
in the definite and a pair of straight lines in the indefinite
case. Summing up:

\bs\noi
{\bf Proposition 1:} {\sl In a sufficiently small neighborhood of a regular
rank-1-bifurcation point (0,0,0) in the $(q,p,\eps)$-space,
the fixed point set consists of a
single point if the bifurcation is definite and it is the union
$A\cup B$ of two smooth 1-dimensional submanifolds, intersecting
at the bifurcation point, if the bifurcation is
indefinite that is in the case of a cross-bifurcation.
Moreover, if $a$ and $b$ denote the tangents to $A$
and $B$ at the bifurcation point, then $a\ne b$, and the
intersection of the plane
spanned by $a$ and $b$ in the $(q,p,\eps)$-space with the plane
defined by fixing the bifurcation parameter value that is with
the $(q,p)$-plane $\eps=0$ is the
1-dimensional eigenspace at the bifurcating fixed point.
The bifurcation is transcritical, if and only if this eigenspace is
different from $a$ and from $b$.}\qed

\bs\noi
Recall that in adapted coordinates the eigenspace is the
$q$-axis and the plane spanned by $a$ and $b$ is $T_0X$, the
$(q,\eps)$-plane.

Now let for a moment $(Q(p,q,\eps), P(p,q,\eps))$ denote any
local symplectic family with a fixed point at $(0,0,0)$, without
assuming it to be rank-1-bifurcating or bifurcating at all. If
there is a differentiable map $\alpha:I\to\IR^2$, defined on an
open interval $I\subset\IR$ around $0$, written in coordinates
as $\alpha(\eps)=(q(\eps),p(\eps))$, such that $\alpha(0)=(0,0)$
and
\begin{equation}
Q(q(\eps),p(\eps),\eps)=q(\eps)\quad\text{and}\quad
P(q(\eps),p(\eps),\eps)=p(\eps)
\end{equation}
for all $\eps\in I$, then we call its graph
\begin{equation}
A:=\{(q(\eps),p(\eps),\eps)\mid\eps\in I\}\subset\IR^3
\end{equation}
a {\nbf fixed point branch} of $(0,0,0)$, regardless of what
other fixed points of the symplectic family might exist.
For instance, if a transcritical bifurcation is considered in a
suitable small neighborhood, then its fixed point set will
consist of two branches $A$ and $B$ as in Proposition 1, but
if an indefinite rank-1-bifurcation is {\it not} trans\-critical, then
only one of its fixed point lines, say $A$, will be a branch,
while $B$, being tangent to $\eps=0$, cannot be
parametrized differentiably by $\eps$.

If again $A$ is a fixed point branch of $(0,0,0)$ in some symplectic family,
then the trace of the Jacobian matrix, taken at each point of
the branch, defines a differentiable real valued function
$\Tr_A:I\to\IR$, so
\begin{equation}
{\rm Tr}_A(\eps):=
Q_q(q(\eps),p(\eps),\eps)+P_p(q(\eps),p(\eps),\eps).
\end{equation}
Then the eigenvalue at $(0,0,0)$ is 1 if and only if
${\rm Tr}_A(0)=2$. Under what additional condition will
$(0,0,0)$ be an indefinite rank-1-bifurcation point, i.e. a
cross-bifurcation?

\bs\noi
{\bf Proposition 2:} {\sl If a fixed point branch $A$ of $(0,0,0)$ in a
symplectic family satisfies ${\rm Tr}_A(0)=2$, then
$(0,0,0)$ will be a cross-bifurcation point if and
only if ${\rm Tr}_A'(0)\ne0$.}

\bs\noi
{\caps Proof:} So let a fixed point branch $A$ of $(0,0,0)$ in a
symplectic family $(Q,P)$ be given, and ${\rm Tr}_A(0)=2$ be
assumed. We have to show (a): If ${\rm Tr}'_A(0)\ne0$, then
$(0,0,0)$ is a cross-bifurcation point, and conversely (b):
If $(0,0,0)$ is a cross-bifurcation point, then
${\rm Tr}'_A(0)\ne0$.

\bs\noi
Proof of (a): We first choose the symplectic $(q,p)$-coordinates
in such a way, that the $q$-axis is contained in the eigenspace,
and the $\eps$-axis $q=p=0$ is tangent to the branch at $(0,0,0)$.
Then we already have
\begin{equation}
\left(
\begin{array}{ccc}
 Q_{ q}\mathord-1& Q_{
p}& Q_\eps\\
 P_{ q}& P_{
p}\mathord-1& P_\eps\\
\end{array}
\right)
=
\left(
\begin{array}{ccc}
0&Q_p&0\\
0&0&0\\
\end{array}
\right)\,,
\end{equation}
at $(0,0,0)$, but we do not yet know if $Q_p\ne0$, and we know nothing about
the Hessian of $P$, except $P_{\eps\eps}(0,0,0)=0$, which
follows from differentiating $P(q(\eps),p(\eps),\eps)=p(\eps)$
twice at $\eps=0$, note $q'(0)=p'(0)=0$ by our choice of
coordinates. So it remains to show that $Q_p\ne0$ and
$P_{q\eps}\ne0$ at the bifurcation point. Now since $(Q,P)$ is a
symplectic family, we have
\begin{equation}
\det\left(
\begin{array}{ll}
Q_q&Q_p\\
P_q&P_p\\
\end{array}
\right)
=
1,
\end{equation}
everywhere, not only at fixed points. So if $u$ denotes any of
our three variables $q$, $p$ or $\eps$, we obtain
\begin{equation}
\det\left(
\begin{array}{ll}
Q_{qu}&Q_p\\
P_{qu}&P_p\\
\end{array}
\right)
+
\det\left(
\begin{array}{ll}
Q_q&Q_{pu}\\
P_q&P_{pu}\\
\end{array}
\right)
=
0
\end{equation}
everywhere, and thus at the special point $(0,0,0)$ we get
\begin{equation}
\det\left(
\begin{array}{ll}
Q_{qu}&Q_p\\
P_{qu}&1\\
\end{array}
\right)
+
\det\left(
\begin{array}{ll}
1&Q_{pu}\\
0&P_{pu}\\
\end{array}
\right)
=
0.
\end{equation}
For later reference, let us call this the

\bs\noi
{\bf Determinant derivative formula:} {\sl If
$Q_q=P_p=1$ and $P_q=0$ at some point in a
symplectic family , then if $u$
is any of the three variables $q$, $p$, or $\eps$ we have
\begin{equation}
Q_{qu}+P_{pu}=Q_pP_{qu}
\end{equation}
at this particular point.}

\bs\noi
Presently we apply the formula for $u:=\eps$. Since
${\rm Tr}_A'(0)=Q_{q\eps}(0,0,0)+P_{p\eps}(0,0,0)$
because of $q'(0)=p'(0)=0$, we have
${\rm Tr}_A'(0)=Q_pP_{q\eps}$ at $(0,0,0)$,
so $Q_pP_{q\eps}\ne0$ by our assumption ${\rm Tr}_A'(0)\ne0$,
and (a) is proved.

\bs\noi
Proof of (b): Now we assume $P_{qq}P_{\eps\eps}-P_{q\eps}^2<0$
at $(0,0,0)$ in adapted coordinates, and we have to show that this implies
${\rm Tr}_A'(0)\ne0$. The branch $A$ must be tangent to the
$(q,\eps)$-plane at the bifurcation point, which means
$p'(0)=0$, and a change to adapted coordinates
of the type
\begin{equation}
\begin{array}{lll}
\widetilde q&=&q-c_1\eps\\
\widetilde p&=&p\\
\end{array}
\end{equation}
with a suitable constant $c_1$ will even make the $\eps$-axis
tangent to $A$, with the Hessian still being indefinite in the
new coordinates. The trace function ${\rm Tr}_A(\eps)$ remains
unchanged anyway, as it is independent of the coordinate choice.
So we may assume $q'(0)=p'(0)=0$ for the branch from the start.
As before, this implies
${\rm Tr}_A'(0)=Q_{q\eps}+P_{p\eps}$
and $P_{\eps\eps}=0$ at $(0,0,0)$. Since we have a
rank-1-bifurcation, we know $Q_p\ne0$, and since it is regular,
$P_{q\eps}\ne0$ follows from $P_{\eps\eps}=0$ at $(0,0,0)$.
The determinant derivative formula is
applicable and shows ${\rm Tr}_A'(0)=Q_pP_{q\eps}$ and hence
${\rm Tr}_A'(0)\ne0$, which completes the proof of (b) and of
Proposition~2.\qed

\bs\noi
As we have seen in Proposition 1, locally at a
{\it transcritical} bifurcation point the fixed point set
consists of two branches $A$ and $B$. Since the functions
${\rm Tr}_A(\eps)-2$ and ${\rm Tr}_B(\eps)-2$ both change sign
at $\eps=0$, both branches change their stability properties from
elliptic to hyperbolic or vice versa, and in fact in opposite
directions. More precisely

\bs\noi
{\bf Proposition 3:} {\sl If $A$ and $B$ are the branches of a
transcritical bifurcation point, then
${\rm Tr}_A'(0)+{\rm Tr}_B'(0)=0$.}

\bs\noi
{\caps Proof:} In adapted coordinates we have
$p_A'(0)=p_B'(0)=0$ anyway, since locally the fixed point set is
contained in the surface $X$ tangent to $p=0$. Changing the
adapted coordinates by a suitable transformation of the type
(21), we can also obtain $q_A'(0)+q_B'(0)=0$. In these
coordinates, we have $P_{q\eps}(0,0,0)=0$.
On the other hand,
${\rm Tr}_A'(0)+{\rm Tr}_B'(0)=
Q_{qq}q_A'(0)+Q_{q\eps}+P_{pq}q_A'(0)+P_{p\eps}
+
Q_{qq}q_B'(0)+Q_{q\eps}+P_{pq}q_B'(0)+P_{p\eps}
=2(Q_{q\eps}+P_{p\eps})$, but this is
$2Q_pP_{q\eps}$
by the determinant derivative formula, and the proposition
follows.\qed

\vskip1.cm
\noi
{\sefont 3 Fork-like bifurcations}

\ms\noi
Let us now have a look at those cross-bifurcations,
which are {\it not} transcritical. Then only one of the two
fixed point lines $A$ and $B$ of Proposition~1 will be a {\it branch}
in the technical sense (14), while the other, say $B$, is
tangent to the $q$-axis. Thus the function $B\to\IR$ given by the
restriction $\eps|B$ of the family parameter to $B$ is critical
at the bifurcation point. It may in fact happen that
$\eps|B\equiv0$, the most degenerate possibility. But here
we introduce a terminology
for the {\it least} degenerate case.

\bs\noi
{\bf Definition 3:} Let
a cross-bifurcation which is not transcritical
be given. Then if the
restriction $\eps|B$ of the bifurcation parameter to one of the
local fixed point lines has a nondegenerate extremum at the
bifurcation point, we speak of a {\nbf fork-like}
bifurcation.\qed

\bs\noi
It is the same condition on $B$ as is required of the fixed
point set in the definition of an
ordinary saddle-node rank-2-bifurcation. But in our case, the
additional fixed point branch $A$ is present, while in a
saddle-node bifurcation the fixed point set is locally just a
single line.

Since $B$ in the non-transcritical case
is tangent to the eigenspace, it can be para\-metrized
as $(q,p_B(q),\eps_B(q))$ by $q$ in adapted coordinates, with
${p_B'(0)=\eps_B'(0)=0}$. The fork-like condition then just
says $\eps_B''(0)\ne0$. In the surface $X$ of (8), the line $B$
locally looks like a parabola up to higher order terms,
tangent to $\eps=0$.

The name is meant to refer to the well-known
{\it pitchfork bifurcations} in symplectic families. These are
bifurcations of fixed points that have eigenvalue $-1$ and as
such do not fall under our heading. But if we consider the
{\it iterated} family instead,
$(Q(Q(q,p,\eps),P(q,p,\eps),\eps),
P(Q(q,p,\eps),P(q,p,\eps),\eps),\eps)$ in coordinates, then the
fixed point gets the eigenvalue $+1$ and the iterated family
bifurcates it fork-like. But being iterations, these are very special
fork-like bifurcations. For instance does the original family
provide a sort of symmetry, which a fork-like bifurcation
in general will not have.

The trace along the fixed point line $B$ of a fork-like
bifurcation is of course well-defined as a
function ${\rm Tr}_B$ on $B$ itself, but not like ${\rm Tr}_A$
as a function of $\eps$, since for small $\eps\ne0$ we find
either none or two points on $B$, and these two may in fact have
different traces. But in an natural way, ${\rm Tr}_B'(0)$ is
still defined and in a fixed relation to ${\rm Tr}_A'(0)$, as
the next proposition shows.

\bs\noi
{\bf Proposition 4:} {\sl  Let $A$ be the fixed point branch of a
fork-like bifurcation and $B$ the fixed point line tangent to
$\eps=0$. Then there exists the limit
\begin{equation}
\lim\limits_{q\to0}\frac{{\rm Tr}_B(q,p_B(q),\eps_B(q))-2}{\eps_B(q)}
=:{\rm Tr}_B'(0)
\end{equation}
and it satisfies
\begin{equation}
{\rm Tr}_A'(0)+{\textstyle\frac12}{\rm Tr}_B'(0)=0\,.
\end{equation}
}

\bs\noi
{\caps Proof:} To determine the limit, we will use second order
Taylor expansion of $Q_q(q,p,\eps)+P_p(q,p,\eps)$ in adapted
coordinates (7). All partial derivatives now taken at $(0,0,0)$.
By the determinant derivative formula we know
\begin{equation}
\begin{array}{lll}
(Q_q+P_p)_q&=&Q_pP_{qq}\\
(Q_q+P_p)_p&=&Q_pP_{qp}\\
(Q_q+P_p)_\eps&=&Q_pP_{q\eps}\,,\\
\end{array}
\end{equation}
and $P_{qq}=0$ since the bifurcation is not transcritical.
Moreover,
\begin{equation}
{\rm Tr}_A'(0)=(Q_q+P_p)_qq_A'(0)+(Q_q+P_p)_pp_A'(0)+
(Q_q+P_p)_\eps
\end{equation}
and $p_A'(0)=0$, so from (24) and $P_{qq}=0$ we also get
\begin{equation}
{\rm Tr}_A'(0)=Q_pP_{q\eps}\,.
\end{equation}
Thus we know the linear terms, and the quotient
becomes
\begin{equation}
\begin{array}{ll}
\displaystyle\frac{{\rm Tr}_B(q,p_B(q),\eps_B(q))-2}{\eps_B(q)}
=&
Q_pP_{qp}\displaystyle\frac{p_B(q)}{\eps_B(q)}+{\rm Tr}_A'(0)\\
\noalign{\ms}
&+
\displaystyle\frac1{\eps_B(q)}\cdot (\text{higher order terms})
\end{array}
\end{equation}
But not many of the higher order terms will contribute to the limit,
since
\begin{equation}
\lim\limits_{q\to0}\frac{p_B(q)}{\eps_B(q)}=
\frac{p_B''(0)}{\eps_B''(0)}\quad\text{and}\quad
\lim\limits_{q\to0}\frac{q^2}{\eps_B(q)}=
\frac2{\eps_B''(0)}\,,
\end{equation}
and as an intermediate result we obtain that the limit
${\rm Tr}_B'(0)$ exists and that
\begin{equation}
{\rm Tr}_B'(0)={\rm Tr}_A'(0)+
\frac{Q_pP_{qp}p_B''(0)+Q_{qqq}+P_{pqq}}{\eps_B''(0)}\,.
\end{equation}
To prove (23), we have to prove that the quotient on the right hand side of
this last equation (29) equals $-3{\rm Tr}_A'(0)$, and thus by
(26) what remains to be shown is
\begin{equation}
3Q_pP_{q\eps}\eps_B''(0)+Q_pP_{qp}p_B''(0)+Q_{qqq}+P_{pqq}=0\,.
\end{equation}
Now we use the fixed point property
\begin{equation}
\begin{array}{l}
Q(q,p_B(q),\eps_B(q))=q\\
P(q,p_B(q),\eps_B(q))=p_B(q)
\end{array}
\end{equation}
of $B$. Differentiating the first equation twice and the second
three times at $q=0$, we get
\begin{equation}
Q_{qq}+Q_p\, p_B''(0)=0
\end{equation}
and
\begin{equation}
P_{qqq}+3P_{q\eps}\,\eps_B''(0)+3P_{qp}\, p_B''(0)=0.
\end{equation}
Inserting this into equation (30), which we have to prove,
we find that (30) is equivalent to
\begin{equation}
Q_{qqq}+P_{pqq}=Q_pP_{qqq}-2Q_{qq}P_{pq}
\end{equation}
at $(0,0,0)$. But this turns out to be a consequence of a second
order determinant derivative formula: we start from the
symplectic property (17), apply $\de^2/\de q^2$,
put in what we know about the first and second order partial
derivatives of $Q$ and $P$ at $(0,0,0)$, and out comes (34).
Thus Proposition~4 is proved.\qed

\bs\noi
From the calculations of the proof, let us preserve the formula
\begin{equation}
\eps_B''(0)=\frac{3Q_{qq}P_{qp}-Q_pP_{qqq}}{3Q_pP_{q\eps}}\,,
\end{equation}
which is a consequence of (32) and (33). Since in its
derivation the non-vanishing of $\eps_B''(0)$ was not used, we
have as a

\bs\noi
{\bf Corollary:} {\sl A non-transcritical
cross-bifurcation is fork-like if and only~if
\begin{equation}
3Q_{qq}P_{qp}\ne Q_pP_{qqq}
\end{equation}
at $(0,0,0)$ in adapted coordinates.}\qed

\vskip1.cm
\noi
{\sefont 4 Fixed point branches given by librating orbits}

\ms\noi
Now let $H:M\to\IR$ be an autonomous Hamiltonian on a
4-dimensional symplectic manifold and $\gamma:\IR\to M$ a
periodic orbit with period $T>0$ of the Hamiltonian flow.
Choose any 3-dimensional
submanifold $\Sigma\subset M$ which is being intersected
transversally by $\gamma$ at time $t=0$. Then on a sufficiently
small neighborhood $U\subset\Sigma$ of $\gamma(0)$ in $\Sigma$,
the Poincar\'e map ${\rm Poinc}:U\to\Sigma$ is well-defined by
following the orbits until they hit $\Sigma$ again after
travelling approximately the period time $T$ of $\gamma$.
Let $E_0:=H(\gamma(0))$ be the energy of the fixed point
$\gamma(0)$ and write $U_E:=U\cap H\inv(E)$ and
$\Sigma_E:=\Sigma\cap H\inv(E)$. After making $\Sigma$ smaller
if necessary, the individual $U_E$ and $\Sigma_E$ will be
nondegenerate subsurfaces of $M$, with the Poincar\'e map defining
a symplectic map $U_E\to\Sigma_E$ for each $E$.
Now we introduce a coordinate $\eps$ for the energy, like
$\eps=E-E_0$, and extend it to local coordinates $(q,p,\eps)$ for
$\Sigma$ such that $(q,p)$ are symplectic coordinates on each $\Sigma_E$.
Then the Poincar\'e map will be a symplectic family, described by
two functions $Q(q,p,\eps)$ and $P(q,p,\eps)$ in these
coordinates. All this is of course well-known, recalled here
only to introduce notation.

Often used is the following sort of `automatic' choice of
$\Sigma$ and the coordinates $q$ and $p$. Let $H=H(x,y,p_x,p_y)$
be given on $M=\IR^4$. If the periodic orbit $\gamma$ satisfies
$\dot y(0)\ne0$, then the 3-dimensional subspace $y=y_0:=y(0)$
in $\IR^4$ can be taken for a start to find $\Sigma$. Of course,
far away from the point $\gamma(0)$ this space may have bad
properties with respect to the Hamiltonian flow. But if we
choose $\Sigma$ as a sufficiently small open neighborhood of $\gamma(0)$
in this 3-space, then not only the Poincar\'e map and the energy
surfaces $\Sigma_E$ will be defined as described, but also
the projection of each $\Sigma_E$ to the $(x,p_x)$-plane will be
symplectic and injective, which means that $x$ and $p_x$ can be
used as the symplectic coordinates $q$ and $p$ on $\Sigma_E$.

The transcritical bifurcations studied in \cite{BT} are related to
straight-line librating orbits of Hamiltonians of the form
\begin{equation}
H(x,y,p_x,p_y)=\textstyle\frac12p_x^2+\frac12p_y^2+V(x,y).
\end{equation}
Let $\gamma$ be such a straight-line periodic orbit, without
loss of generality projecting to the $y$-axis, oscillating there
between two values $y_1<y_2$, with $x(0)=0$ and $y_1<y(0)<y_2$
and an energy $E_0$. As just recalled, we use
\begin{equation}
\begin{array}{l}
q:=x\\
p:=p_x\\
\eps:=E-E_0
\end{array}
\end{equation}
as coordinates to describe the `automatic' Poincar\'e map by a
symplectic family $(Q(q,p,\eps),P(q,p,\eps))$, with the fixed
point $(0,0,0)$ corresponding to the orbit $\gamma$. What do we
know about this symplectic family?

Straight-line librating orbits in Hamiltonian systems of type
(37) always come in one-parameter families. The given orbit
$\gamma$ satisfies
\begin{equation}
\dot p_x(t)=-\frac{\de V}{\de x}(0,y(t))\equiv0,
\end{equation}
hence at least we know $\frac{\de V}{\de x}(0,y)=0$ for all $y$ with
$y_1\le y\le y_2$. If $V$ is polynomial or real-analytic, this implies
\begin{equation}
\frac{\de V}{\de x}(0,y)=0\quad\text{for all}\quad y\in\IR\,.
\end{equation}
In this case the whole $(y,p_y)$-plane $x=p_x=0$ in $\IR^4$ is
invariant under the Hamiltonian flow, and the flow lines are
just those of the
1-dimensional Hamiltonian system on the $(y,p_y)$-plane given by
\begin{equation}
h(y,p_y)=\textstyle\frac12p_y^2+V(0,y).
\end{equation}
By standard regularity arguments, the orbit $\gamma$ must be
embedded in a family of neighboring closed orbits of this
1-dimensional system, one for each energy in an interval
$(E_1,E_2)$ with $E_1<E_0<E_2$. Thus we have a family of
librating orbits of the original 2-dimensional system on the
$y$-axis over\footnote{Here we assumed (40). What if $\frac{\de V}{\de x}(0,y)$ may be
non-zero for $y>y_2$ or $y<y_1$? Then the neighboring closed
orbits on the {\it inside} of $\gamma$, those with energies
$E_1<E<E_0$, still are librating orbits of the original system as well,
and define a fixed point branch in the Poincar\'e family with
$q=p=0$ and $\eps_1<\eps<0$. The orbits of the 1-dimensional
system which are on the outside of $\gamma$, with energies
$E_0<E<E_2$, need not be orbits of the 2-dimensional system. But
then we simply shift attention from $E_0$ to an energy $E_0'$
between $E_1$ and $E_0$. So in any case, we may assume a fixed
point branch $A$ as in (42) in the Poincar\'e symplectic family be
given, describing a family of neighboring straight-line
librating orbits on the $y$-axis.}
$(E_1,E_2)$. In our Poincar\'e symplectic family it
corresponds to a fixed point branch
\begin{equation}
A:=\{(0,0,\eps)\mid\eps_1<\eps<\eps_2\}.
\end{equation}
For a branch $A$, the trace function
${\rm Tr}_A:(\eps_1,\eps_2)\to\IR$ is defined, and if it happens
to be that ${\rm Tr}_A(0)=2$ and ${\rm Tr}_A'(0)\ne0$, then by
Proposition~2 we have a cross-bifurcation point, which
may or may not be transcritical. Examples of both cases have
been presented numerically in \cite{BT}.

Hamiltonians of type (37) have time-reversal symmetry, and in
\cite{BT} the destruction of a transcritical bifurcation by a
symmetry breaking perturbation is described. What role does
time-reversal symmetry play in transcritical bifurcations?
Could we {\it preserve} a transcritical bifurcation under a
symmetry breaking perturbation, and conversely, can a
transcritical bifurcation be {\it destroyed} by a symmetric
perturbation? The answer to the first question is yes, as we are
going to show now.

\vskip1.cm
\noi
{\sefont 5 Transcritical bifurcation in unsymmetric Hamiltonian
systems}

\ms\noi
Again we start from a Hamiltonian
$H(x,y,p_x,p_y)=\frac12p_x^2+\frac12p_y^2+V(x,y)$ with a
straight-line librating orbit $\gamma$, embedded in a family of
such orbits, all projecting to the $y$-axis and thus constituting a
fixed point branch $A=0\times0\times(\eps_1,\eps_2)$ in the
Poincar\'e symplectic family as described above, with $(0,0,0)$
representing $\gamma$. For convenience and with little loss of
generality we assume (40), that was
$\frac{\de V}{\de x}(0,y)=0$ for {\it all} $y\in\IR$. But our
main assumption now will be that $(0,0,0)$ is a
cross-bifurcation point in the Poincar\'e family. We then speak of
$\gamma$ as of a {\nbf cross-libration} or of a {\nbf
transcritical libration}, if this cross-bifurcation happens to
be transcritical.

Now let the Hamiltonian depend on an additional small parameter
$\delta$. Such a function $\overline H(x,y,p_x,p_y,\delta)$
will be called a {\nbf perturbation} of $H$ or of the
transcritical libration, if
$\overline H(x,y,p_x,p_y,0)=H(x,y,p_x,p_y)$.
Then we have a two-parameter Poincar\'e symplectic family
\begin{equation}
(\overline Q,\overline P)=
(\overline Q(q,p,\eps,\delta),
\overline P(q,p,\eps,\delta)),
\end{equation}
defined on an open neighborhood of $(0,0,0,0)$ in the
$(q,p,\eps,\delta)$-space $\IR^4$ in the usual way, with
$q$ and $p$ from $x$ and $p_x$.

\bs\noi
{\bf Definition 4:} We say that a perturbation $\overline H$ of
a cross-libration ist {\nbf cross-preser\-ving}, if
there are differentiable functions $q(\delta)$, $p(\delta)$ and
$\eps(\delta)$, defined on an interval $(-\delta_0,\delta_0)$,
with $q(0)=p(0)=\eps(0)=0$, such that for any fixed
$\delta\in(-\delta_0,\delta_0)$ the point
$(q(\delta),p(\delta),\eps(\delta))$ is a
cross-bifurcation point in the symplectic $\eps$-family given by
$(\overline Q,\overline P)$ at the fixed $\delta$.\qed

\bs\noi
Note that a transcritical libration will stay transcritical
under a cross-preserving deformation, for small enough
$\delta_0$.

\bs\noi
There is a simple strategy to find cross-preserving
perturbations. All we have to do is to make sure that for all
sufficiently small $\delta$ the $(y,p_y)$-plane
$0\times\IR\times0\times\IR$ is invariant under
the Hamiltonian flow on $\IR^4$ that is
defined by $\overline H$ for fixed $\delta$. In other words, if
\begin{equation}
\frac{\de\overline H}{\de p_x}(0,y,0,p_y,\delta)=0
\quad\text{and}\quad
\frac{\de\overline H}{\de x}(0,y,0,p_y,\delta)=0
\end{equation}
for sufficiently small $\delta$ and all $(y,p_y)$, then
$\overline H$ will preserve the transcriticality. Why?
Because then the closed orbits neighboring $\gamma$ will define
a 2-parameter fixed point branch $\overline A$ of
$(\overline Q,\overline P)$ given by
functions $q(\eps,\delta)$ and $p(\eps,\delta)$, defined on a
neighborhood of $(0,0)$ in the $(\eps,\delta)$-plane. We know
\begin{equation}
{\rm Tr}_{\overline A}(0,0)=2
\quad\text{and}\quad
\frac{\de{\rm Tr}_{\overline A}}{\de\eps}(0,0)\ne0
\end{equation}
by assumption, thus by the implicit function theorem, we get a
differentiable function $\eps(\delta)$ with $\eps(0)=0$ and
\begin{equation}
{\rm Tr}_{\overline A}(\eps(\delta),\delta)=2
\quad\text{and}\quad
\frac{\de{\rm Tr}_{\overline A}}{\de\eps}(\eps(\delta),\delta)\ne0
\end{equation}
for sufficiently small $\delta$. By Proposition~2, each
$(q(\eps(\delta),\delta),p(\eps(\delta),\delta),\eps(\delta))$
is a cross-bifurcation point in the symplectic family of
the corresponding fixed $\delta$ (and if it is transcritical
for $\delta=0$, it must remain transcritical for sufficiently
small $\delta$). So the strategy is sound, and we use it to
derive the following proposition.

\bs\noi
{\bf Proposition 5:} {\sl Let
$H(x,y,p_x,p_y)=\frac12p_x^2+\frac12p_y^2+V(x,y)$ be a
Hamiltonian with a cross-libration on the $y$-axis
and ${\de V}/{\de x}(0,y)=0$ for all $y\in\IR$. Let
$F(x,y,p_x,p_y)$ be an arbitrary differentiable function
with
\begin{equation}
\frac{\de F}{\de p_x}(0,y,0,p_y)=0
\quad\text{and}\quad
\frac{\de F}{\de x}(0,y,0,p_y)=0
\end{equation}
Then
\begin{equation}
\overline H(x,y,p_x,p_y,\delta):=H(x,y,p_x,p_y)+
\delta\mkern4mu F(x,y,p_x,p_y)
\end{equation}
is a cross-preserving perturbation, since it
obviously satisfies} (44).\qed

\bs\noi
It is a matter of definition, which Hamiltonians should be
counted as `unsymmetric' in this context. Making $H$ unsymmetric
by a bump far away from the transcritical libration would not
illuminate the relation between symmetry and transcriticality.
But Proposition~5 shows transcritical bifurcation to occur in the
Hamiltonian system
$H(x,y,p_x,p_y)+\delta\mkern4mu F(x,y,p_x,p_y)$ on $\IR^4$ for fixed
sufficiently small $\delta\ne0$. Since the only condition on
$F(x,y,p_x,p_y)$ is the vanishing of the partial derivatives by
$x$ and $p_x$ along the $(y,p_y)$-plane,
we
may fairly say yes, transcritical bifurcation can occur in
{\it unsymmetric} autonomous 2-dimensional Hamiltonian systems.

\vskip1.cm
\noi
{\sefont 6 Destruction of cross-bifurcations}

\ms\noi
Once more we start with a cross- or a transcritical
libration $\gamma$ on the
$y$-axis of a Hamiltonian system
$H(x,y,p_x,p_y)=\frac12p_x^2+\frac12p_y^2+V(x,y)$
with $\frac{\de V}{\de x}(0,y)=0$ for all $y\in\IR$. This time
we ask by which perturbations of $H$ the cross-property
can be destroyed. By destruction we mean more than
just non-preservation:

\bs\noi
{\bf Definition 5:} We say that a perturbation $\overline H$ of
a cross-libration is {\nbf cross-destroy\-ing}, if
there is a neighborhood of $(0,0,0,0)$ in
$(q,p,\eps,\delta)$-space, in which $(0,0,0,0)$ is the only
cross $\eps$-bifurcation point of $(\overline Q,\overline P)$.
\qed

\bs\noi
We will first prove a general destruction criterion for
cross-bifurcation in symplectic families.

\bs\noi
{\bf Proposition 6:} {\sl Let $(Q(p,q,\eps),P(p,q,\eps))$ be a
symplectic family with a cross-bifurcation at $(0,0,0)$.
Extend it to a symplectic family
$(\overline Q(p,q,\eps,\delta),\overline P(p,q,\eps,\delta))$
depending on a second parameter $\delta$, with
$\overline Q(p,q,\eps,0)=Q(p,q,\eps)$ and
$\overline P(p,q,\eps,0)=P(p,q,\eps)$. Then if the vector
$(\overline Q_\delta(0,0,0,0),\overline P_\delta(0,0,0,0))$ does
not belong to the eigenspace of $(Q,P)$ at the bifurcation
point, there is a neighborhood of $(0,0,0,0)$ in the
$(q,p,\eps,\delta)$-space, in which it is the only
$\eps$-cross-bifurcation point.}

\bs\noi
{\caps Proof:} The assumption simply means that
at $(0,0,0,0)$ the Jacobian
matrix
\begin{equation}
\left(
\begin{array}{cccc}
Q_q\mathord-1&Q_p&Q_\eps&\overline Q_\delta\\
P_q&P_p\mathord-1&P_\eps&\overline P_\delta\\
\end{array}
\right)
\end{equation}
of $(\overline Q-q,\overline P-p)$
is of rank 2, as can be seen most easily in adapted coordinates,
where the condition reduces to $\overline P_\delta\ne0$.
So the fixed point set $\overline F$, locally at $(0,0,0,0)$, is
a smooth surface, and since it contains the two fixed point
lines $A$ and
$B$ of the given cross-bifurcation,
its tangent plane at the point must  be spanned by $a$ and
$b$ und thus be the $(q,\eps)$-plane in adapted coordinates.
In particular, the restriction $\delta|\overline F$ is singular
at $(0,0,0,0)$.

To prove the proposition, it would be sufficient to show that
this is an {\it isolated} singularity of $\delta|\overline F$.
Because then at any other fixed point
$(q_0,p_0,\eps_0,\delta_0)$ in a neighborhood, the set of fixed
points with $\delta=\delta_0$ will be locally a smooth
1-dimensional submanifold of $\overline F$ and hence the point
can't be an  $\eps$-cross-bifurcation point.

To find out, we calculate the matrix of the Hessian quadratic
form of $\delta|\overline F$ at the singularity. Straightforward
calculation using Lagrangian multipliers gives
\begin{equation}
{\rm Hess}_{\delta|\overline F}(0,0,0,0)=
-\frac1{\overline P_\delta}\left(
\begin{array}{ll}
P_{qq}&P_{q\eps}\\
P_{q\eps}&P_{\eps\eps}\\
\end{array}
\right)
\end{equation}
for this matrix, which is now seen to be
nondegenerate by assumption, hence the
singularity is isolated and Proposition~6 is proved.\qed

\bs
\noi
Let $\gamma$ be a cross-libration on the
$y$-axis of a Hamiltonian system
\begin{equation}
H(x,y,p_x,p_y)=\frac12p_x^2+\frac12p_y^2+V_0(x,y)
\end{equation}
with $\frac{\de V_0}{\de x}(0,y)=0$ for all $y\in\IR$.
For an attempt to apply Proposition~6 to perturbations
\begin{equation}
\overline H(x,y,p_x,p_y,\delta):=H(x,y,p_x,p_y)+
\delta\mkern4mu F(x,y,p_x,p_y)
\end{equation}
with an arbitrary perturbation term $F(x,y,p_x,p_y)$,
we will have to understand how
$(\overline Q_\delta(0,0,0,0),\overline P_\delta(0,0,0,0))$ does
depend on $F(x,y,p_x,p_y)$.

\bs\noi
{\bf Proposition 7:} {\sl Let $T>0$ denote the period of the
cross-libration $\gamma$ and $y(t)$ its $y$-component.
Define the function $f(y)$ by
\begin{equation}
f(y):=\frac{\de^2V_0}{\de x^2}(0,y).
\end{equation}
and let $g_1(y,p_y)$ and $g_2(y,p_y)$ be given by
\begin{equation}
\begin{array}{l}
g_1(y,p_y):=\phantom{-}\displaystyle\frac{\de F}{\de
p_x}(0,y,0,p_y),\\
\noalign{\ms}
g_2(y,p_y):=-\displaystyle\frac{\de F}{\de x}(0,y,0,p_y).
\end{array}
\end{equation}
Let $(\xi(t),\eta(t))$ be the solution of
\begin{equation}
\begin{array}{lrll}
\dot\xi&-\eta&=&g_1(y(t),\dot y(t))\\
\dot\eta&+f(y(t))\xi&=&g_2(y(t),\dot y(t))
\end{array}
\end{equation}
with the initial condition $\xi(0)=\eta(0)=0$. Then
$\overline H$ satisfies the destruction criterion of
Proposition~6 if and only if $(\xi(T),\eta(T))$ is
not contained in the eigenspace of the Poincar\'e map 
of the undisturbed system at the bifurcation point.}\qed

\bs\noi
What would it take to apply Proposition~7 numerically? First of
all, one has to know the component $y(t)$ and its derivative
$\dot y(t)=p_y(t)$ of the cross-librating
orbit~$\gamma$. Secondly, we will need the fundamental system
$(\vf(t),\psi(t))$ of the homogeneous equation
\begin{equation}
\ddot\xi+f(y(t))\xi=0
\end{equation}
with initial conditions $\vf(0)=\dot\psi(0)=1$ and
$\dot\vf(0)=\psi(0)=0$,
and finally we would have to calculate the two numbers
\begin{equation}
\begin{array}{l}
c_1(T):=\phantom{-}\int\limits_0^T(\dot\psi(t)g_1(y(t),\dot y(t))-
\psi(t)g_2(y(t),\dot y(t))\,dt\\
c_2(T):=-\int\limits_0^T(\dot\vf(t)g_1(y(t),\dot y(t))-
\vf(t)g_2(y(t),\dot y(t))\,dt.
\end{array}
\end{equation}
This then is all we need, because as can be shown, the Jacobian
matrix of the Poincar\'e map at the bifurcation point turns out to
be
\begin{equation}
\left(
\begin{array}{ll}
Q_q&Q_p\\
P_q&P_p\\
\end{array}
\right)=
\left(
\begin{array}{ll}
\vf(T)&\psi(T)\\
\dot\vf(T)&\dot\psi(T)\\
\end{array}
\right),
\end{equation}
and standard calculation of $(\xi(T),\eta(T))$ leads to the
following

\bs\noi
{\bf Corollary:} {\sl The perturbation term $F(x,y,p_x,p_y)$
satisfies the destruction criterion of
Proposition~6 if and only if
\begin{equation}
\left(
\begin{array}{cc}
\vf(T)-1&\psi(T)\\
\dot\vf(T)&\dot\psi(T)-1\\
\end{array}
\right)
\left(
\begin{array}{c}
c_1(T)\\
c_2(T)
\end{array}
\right)\ne0.
\end{equation}
}\qed

\bs\noi
Application to
$F(x,y,p_x,p_y)=xp_y-yp_x$ leads to the condition
\begin{equation}
\left(
\begin{array}{cc}
\vf(T)-1&\psi(T)\\
\noalign{\bs}
\dot\vf(T)&\dot\psi(T)-1\\
\end{array}
\right)
\left(
\begin{array}{r}
\int\limits_0^Ty(t)\dot\psi(t)\,dt\\
-\int\limits_0^Ty(t)\dot\vf(t)\,dt
\end{array}
\right)\ne0.
\end{equation}

%\bs

For a more systematic approach to the calculation of partial 
derivatives of the Poincar\'e map in connection with straight-line 
librations, see \cite{jaen2}.

\end{document}